\documentclass[13pt,tikz,border=8pt,multi]{amsart}
\usepackage{microtype,amsmath,hyperref,amsthm,tikz,forest,array}
\usetikzlibrary{shapes,arrows,trees,calc}
\usetikzlibrary{shadows}
\newtheorem{theorem}{Theorem}[section]

\newtheorem{original conjecture}[theorem]{Original conjecture}
\theoremstyle{definition}
\newtheorem{definition}[theorem]{Definition}
\newtheorem{example}[theorem]{Example}

\newtheorem{cor}[theorem]{Corollary}
\theoremstyle{remark}
\newtheorem{remark}[theorem]{Remark}
\newtheorem{notation}{Notation}

\numberwithin{equation}{section}

\usepackage[square,numbers]{natbib}

\begin{document}
	\title{Creative proofs in combinations}
	\author{Mohammad Arab}
	\address{Department of Mathematics,Lorestan University, Khoramabad, Iran}
	\email{arab.mohammad20000@gmail.com}
	\curraddr{}
	\date{}
	\dedicatory{}
	\maketitle
	\begin{abstract}
		In this article, we present four issues and provide a creative and concise proof for each of them. The four issues are:
		\begin{enumerate}	
			\item
			Inequality $\frac{1}{\sqrt{n\pi+\frac{\pi}{2}}} < \frac{\binom{2n}{n}}{2^{2n}} < \frac{1}{\sqrt{n\pi}} $
			\item
			A special case of Jonathan Wilde's problem 
			\item
			Combination series
			\item
			A feature of powerful numbers. 
			\end{enumerate}	
	\end{abstract}
\section{Introduction}
In the second section, by proving a theorem, we prove the following equation:
\begin{equation}\label{sara2}
\left(\sum_{k=1}^{n}\frac{(f(k))^2}{2k-1}\right)+(2n+1)(f(n))^2=1 \:\:,\:\:f(k):=\prod_{i=1}^k \left(1-\frac{1}{2i}\right)
\end{equation}
The special case of the Wallis integral that we use to prove equation\eqref{sara2} is as follows\cite{A}:
\[
\int_{0}^{\frac{\pi}{2}} \cos^{2k}(x)\, dx=\prod_{i=1}^k \left(1-\frac{1}{2i}\right) \times\frac{\pi}{2}
\]
And then using the theorem, we infer the following inequality:
\begin{equation}\label{sara1}
\frac{1}{\sqrt{n\pi+\frac{\pi}{2}}} < \frac{\binom{2n}{n}}{2^{2n}} < \frac{1}{\sqrt{n\pi}}
\end{equation}
We use the inequality in Jonathan Wilde's problem.
\subsection*{Jonathan Wilde's problem}
\begin{definition}
$a(n)$ is the number of ways to draw $n$ circles in the affine plane.
\end{definition}
Jonathan Wilde stated problem $a(n)$ by setting conditions:\\
 Two circles must be disjoint or meet in two distinct points (tangential contacts are not permitted), and three circles may not meet at a point.

The sequence was proposed by Jonathan Wild, a professor of music at McGill
University, who found the values $a(1)=1$, $a(2)=3$, $a(3)=14$, $a(4)=173$, and, jointly with Christopher Jones, $a(5)=16951$\cite{R}.\\
It is sequence A250001 in the \textit{On-Line Encyclopedia of Integer Sequences} (OEIS)\cite{C}.\\
We prove a special case of this problem in this article, which we express in the following definition:
\begin{definition}
	$B(n)$ is the number of ways to draw $n$ circles in the affine plane, so that no two circles are neither intersecting nor tangent to each other.  
\end{definition}
We express some of the values of $B(n)$ in sequence: 
\[
B(0)=1,\:\:B(1)=1,\:\:B(2)=2,\:\:B(3)=4,\:B(4)=9,\:\:B(5)=20,...
\]
In the third section, we state two formulas for $B(n)$.\\
In theorems \ref{sara 79} and \ref{sara 77}, we count the number of arrangements of $n$ circles on a plane so that no two circles are neither intersecting nor tangent to each other. The simplest solution we can use to prove the theorems is positive integer partitions.\\
In estimating $B(n)$ in another way, we use inequality\eqref{sara1}.\\

We consider three series from the book \textit{Principles And Techniques In Combinatorics}\cite{T}, which we prove in general in the fourth section and express them in this section.\\
We present the two series in question as follows:
\begin{equation}\label{sara 189}
	\sum_{r=0}^{n} \frac{1}{r+1} \binom{n}{r}=\frac{1}{n+1}\left(2^{n+1}-1\right)
\end{equation}
\begin{equation}\label{sara 189-2}
	\sum_{r=0}^{n} \frac{(-1)^r}{r+1}\binom{n}{r}=\frac{1}{n+1}
\end{equation}
Now we express the third series. Putnam proved the following series in 1962\cite{T}:
\[\sum_{r=1}^{n}r^2\binom{n}{r}=n(n+1)2^{n-2}\]
And  two Chinese teachers, Wei Guozhen and Wang Kai, in 1988, showed that\cite{T}:
\begin{equation}\label{sara 13}
	\sum_{r=1}^{n}r^k\binom{n}{r} =\sum_{i=1}^{k} S(k,i) . P_{i}^{n} . 2^{n-i}
\end{equation}
Where $k\le n$ and $S(k,i)$ are the Stirling numbers of the second kind.\\
We prove the generalization of the series\eqref{sara 13} using derivative formulas.\\
There is a property of sequence $(1^n , 2^n ,3^n ,...)$ that we prove in the fifth section.\\
For each $n$ belonging to natural numbers, if we continue the difference of the sentences $n$ steps, we will reach a fixed sequence of $n!$.
\begin{example}
	Consider the sequence of power numbers for power $3$:
	\begin{table}[ht]
		\caption{$a_{k}=k^{3}$}
		\renewcommand\arraystretch{1.5}
		\noindent\[
		\begin{array}{lccccccccccccr}
			{a_{1}}&{,}&{a_{2}}&{,}&{a_{3}}&{,}&{a_{4}}&{,}&{a_{5}}&{,}&{a_{6}}&{,}&{}&{\cdots}\\
			
			{1}&{,}&{8}&{,}&{27}&{,}&{64}&{,}&{125}&{,}&{216}&{,}&{}&{\cdots}\\
			
			{}&{7}&{,}&{19}&{,}&{37}&{,}&{61}&{,}&{91}&{,}&{127}&{,}&{\cdots}\\
			
			{}&{}&{12}&{,}&{18}&{,}&{24}&{,}&{30}&{,}&{36}&{,}&{}&{\cdots}\\
			
			{}&{}&{}&{6}&{,}&{6}&{,}&{6}&{,}&{6}&{,}&{6}&{,}&{\cdots}\\
			{}&{}&{}&{3!}&{,}&{3!}&{,}&{3!}&{,}&{3!}&{,}&{3!}&{,}&{\cdots}\\
		\end{array}
		\]
	\end{table}
\end{example}
After three steps of difference  of sentences, we reached the fixed sequence $3!$.

\section{Binomial inequality}

\begin{theorem}\label{sara 555}
		\[
	\lim_{n \to \infty}\left(1-\sum_{k=1}^{n}\frac{(f(k))^2}{2k-1}\right)=\lim_{n \to \infty}(2n+1)(f(n))^2=\frac{2}{\pi}
	\]
	We express $f(k)$ as a function in the following form:
	\[
	f(k):=\prod_{i=1}^k \left(1-\frac{1}{2i}\right).
	\]
\end{theorem}
\begin{proof}
	
\textit{The proof has two parts.}\\

In the first part, we prove the following limit:
\[
\lim_{n \to \infty}\left(1-\sum_{k=1}^{n}\frac{(f(k))^2}{2k-1}\right)=\frac{2}{\pi}
\]
And in the second part we prove:
\[
1-\sum_{k=1}^{n}\frac{(f(k))^2}{2k-1}=(2n+1)(f(n))^2.
\]
\textit{Proof of the first part.}\\

The following integral is known as the Wallis integral:
\[
\int_{0}^{\frac{\pi}{2}} \cos^{2k}(x)\, dx=\prod_{i=1}^k \left(1-\frac{1}{2i}\right) \times\frac{\pi}{2}
\]
We rewrite the Wallis integral using the limit definition:
\[
\lim_{b \to 0^{+}}\int_{b}^{\frac{\pi}{2}} \cos^{2k}(x)dx=f(k)\times\frac{\pi}{2}
\]
Now consider the following integral:
\begin{equation} \label{sara 1}
	g(b)=\int_{b}^{\frac{\pi}{2}} \sqrt{1-\cos^{2}(x)}\:\: dx
\end{equation}
A solution for calculating the limit of the above expression is as follows:
\begin{equation}\label{sara 5}
	\lim_{b \to 0^{+}} g(b)=\lim_{b \to 0^{+}}\int_{b}^{\frac{\pi}{2}} \sqrt{\sin^{2}(x)}\:\: dx=\lim_{b \to 0^{+}} \left[-\cos(x)\right]_{b}^{\frac{\pi}{2}}=1
\end{equation}
We rewrite integral \ref{sara 1} and calculate its limit in the following form:
\[
\lim_{b \to 0^{+}} g(b)=\lim_{b \to 0^{+}}\int_{b}^{\frac{\pi}{2}}\left(1-\cos^{2}(x)\right)^{\frac{1}{2}}\, dx \quad,\quad \left(0<b\le x \le \frac{\pi}{2} \right)
\]
According to the set limits we have:
\[
|-\cos^2(x)|<1
\]
Using binomial expansion we can express:
\[
\lim_{b \to 0^{+}}g(b)=\lim_{b \to 0^{+}}\int_{b}^{\frac{\pi}{2}}\left(1+\sum_{k=1}^{\infty}\binom{\frac{1}{2}}{k}\left(-\cos^{2}(x)\right)^{k}\right)\, dx
\]
We rewrite the limit in the following form:
\[
\lim_{b \to 0^{+}} g(b)=\lim_{b \to 0^{+}}\left(\frac{\pi}{2}-b\right)+\left(\sum_{k=1}^{\infty}\binom{\frac{1}{2}}{k}(-1)^{k}\left(\lim_{b \to 0^{+}}\int_{b}^{\frac{\pi}{2}}\cos^{2k}(x)\, dx\right)\right)
\]
We calculate the value of $\binom{\frac{1}{2}}{k}$  separately below:
\[
\binom{\frac{1}{2}}{k}=\frac{(-1)^{k-1}}{2k-1}\times\frac{1\times3\times\cdots\times(2k-1)}{2\times4\times\cdots\times(2k)}=\frac{(-1)^{k-1}}{2k-1}f(k)
\]
Now we place this value in the limit:
\[
\lim_{b \to 0^{+}}g(b)=\left(\frac{\pi}{2}+\sum_{k=1}^{\infty}\frac{(-1)^{2k-1}}{2k-1}\times f(k)\times\left(f(k)\times\frac{\pi}{2}\right)\right)
\]
According to \ref{sara 5} it can be stated that:
\[
1=\frac{\pi}{2}\left(1-\sum_{k=1}^{\infty}\frac{(f(k))^2}{2k-1}\right)
\]
Therefore:
\[
\left(1-\sum_{k=1}^{\infty}\frac{(f(k))^2}{2k-1}\right)=\frac{2}{\pi}.
\]
\\
\textit{Proof of the second part.}\\

We prove this part using the method of mathematical induction, now consider the two sequences as follows:
\[
a_{n}=\left(1-\sum_{k=1}^{n}\frac{(f(k))^2}{2k-1}\right) \quad,\quad b_{n}=(2n+1)(f(n))^2
\]
Using induction we will show that:
\[
a_{n}=\left(1-\sum_{k=1}^{n}\frac{(f(k))^2}{2k-1}\right)=(2n+1)(f(n))^2=b_{n}.
\]
Induction base:
\[
a_{1}=b_{1}=\frac{3}{4}.
\]
Induction assumption:
\[
a_{n}=b_{n}
\]
And we try to get the correctness of induction by considering the assumption of induction correctly:
\[
a_{n+1}=b_{n+1}.
\]
By performing calculations it can be obtained that:
\[
a_{n+1}=a_{n}-\frac{(f(n+1))^{2}}{2n+1}
\]
According to the induction assumption we have:
\[
a_{n+1}=b_{n}-\frac{(f(n+1))^{2}}{2n+1}=(2n+1)(f(n))^{2}-\frac{(f(n+1))^{2}}{2n+1}
\]
Instead of $(f(n+1))^{2}$ we put the following relation: 
\[
(f(n+1))^{2}=\left(\left(\frac{2n+1}{2(n+1)}\right)\times f(n)\right)^{2}
\]
So we will have:
\[
a_{n+1}=(2n+1)(f(n))^{2}-\frac{(2n+1)(f(n))^{2}}{{(2n+2)^{2}}}=(2n+1)(f(n))^{2}\left(1-\frac{1}{(2n+2)^2}\right)
\]
Rewrite the phrase as follows:
\[
a_{n+1}=(2n+1)(f(n))^{2}\left(\frac{(2n+1)(2n+3)}{(2n+2)^2}\right)=\left(\frac{2n+1}{2n+2}f(n)\right)^{2}\times(2n+3)
\]
Therefore:
\[
a_{n+1}=(2n+3)(f(n+1))^{2}=(2(n+1)+1)(f(n+1))^{2}
\]
The sentence of induction is imposed:
\[
a_{n+1}=b_{n+1}.
\]

\textit{Based on the two parts of the proof, we will have:}

\[
\lim_{n \to \infty}a_{n}=\lim_{n \to \infty}b_{n}=\frac{2}{\pi}.
\]	
\end{proof}
\begin{remark}According to the stated theorem, the following equation will be obtained: 
\[
\left(\sum_{k=1}^{n}\frac{(f(k))^2}{2k-1}\right)+(2n+1)(f(n))^2=1.
\]
\end{remark}
\begin{cor}
	\[
	\sum_{k=1}^{\infty} \frac{f(k)(2k)!!}{(2k-1)(2k+1)!!}=\sum_{k=1}^{\infty}\frac{1}{4k^{2}-1}=\frac{\pi}{2}
	\]
	\[
	\sum_{k=1}^{\infty} \frac{f(k)}{2k-1}=1
	\]
\end{cor}
\begin{cor}
The following inequality holds for every $n$ belonging to natural numbers:
\begin{equation}\label{sara 3}
	\frac{1}{\sqrt{n\pi+\frac{\pi}{2}}} < \frac{\binom{2n}{n}}{2^{2n}} < \frac{1}{\sqrt{n\pi}}.
\end{equation}
\end{cor}

\begin{proof}
Consider the following sequence:
\[
a_{n}=(2n+1)(f(n))^2
\]
Sequence $ a_{n} $ is a descending sequence, so according to theorem \ref{sara 555} we can say:
\[
a_{n}=(2n+1)(f(n))^{2}>\frac{2}{\pi} 
\]
So we will have:
\[
(f(n))^{2}>\frac{\frac{2}{\pi}}{2n+1} 
\]
Inequality can be rewritten in another form:
\[
f(n)>\frac{1}{\sqrt{n\pi+\frac{\pi}{2}}}.
\]
To prove the other part of the inequality, we consider the following sequence:
\[
b_{n}=(2n)(f(n))^2
\]
The limit of sequence  can be obtained in the following way:
\[
\lim_{n \to \infty}b_{n}=\lim_{n \to \infty}a_{n}-\lim_{n \to \infty}(f(n))^{2}
\]
We got the limit of the sequence $ a_{n} $ in theorem \ref{sara 555} now to get the limit of the sequence $ (f(n))^{2} $, we get $ (f(n))^{2} $:
\[
(f(n))^{2}=\left(\frac{1\times3\times\cdots\times(2n-1)}{2\times4\times\cdots\times(2n)}\right)^{2}
\]
\[
(f(n))^{2}=\left(\frac{(1\times3\times\cdots\times(2n-1))(2\times4\times\cdots\times(2n))}{2^{2n}(n!)^{2}}\right)^{2}
\]
\[
(f(n))^{2}=\left(\frac{(2n)!}{2^{2n}(n!)^{2}}\right)^{2}=\left(\frac{\binom{2n}{n}}{2^{2n}}\right)^{2}
\]
Therefore, using the Stirling approximation\cite{B}, the limit of $ (f(n))^{2} $ can be obtained as follows:
\[
\lim_{n \to \infty}(f(n))^{2}=\lim_{n \to \infty}\left(\frac{1}{\sqrt{n\pi}2^{2n}}\right)^{2}=0
\]
In this case, the limit of $ b_{n} $ will be obtained:
\[
\lim_{n \to \infty}a_{n}=\lim_{n \to \infty}b_{n}=\frac{2}{\pi}
\]
So we will have:
\[
b_{n}=2n(f(n))^{2} < \frac{2}{\pi}
\]
We rewrite this inequality in another form:
\[
f(n) < \frac{1}{\sqrt{n\pi}}.
\]	
\end{proof}
\section{No. of arrangements of n circles in the plane}
\begin{definition}
	A partition of a positive integer $n$ is a set of positive integers whose sum is $n$. Since the ordering is immaterial, we may regard 
	a partition of $n$ as a finite nonincreasing sequence $ n_1 \ge n_2 \ge \cdots \ge n_l$ of positive integers such that $\sum_{i=1}^l
	n_i= n$. So If $n =n_1 +n_2 +\cdots+n_l$ is a partition of $n$, we say that $n$ is partitioned into $l$ parts of sizes $n_1,n_2, ... ,n_l$ respectively\cite{T}.
\end{definition}
We denote the number of  partitions of $n$ where $n$ is a natural number by $p(n)$.
\begin{definition}
	$p_{1}(n,k) $ is the number of partitions of $n$ whose largest size is $k$.
\end{definition}
\begin{definition}
	The function $p_{2}(n,k) $ is a function of counting the number of states  $n$ of the circle, so that the maximum $k-1$ of the circle is inside at least one circle.
\end{definition}
It is clear from the definition of $p_{2}(n,k) $ and $p_{1}(n,k) $ that for any natural number such as $ n $:
\[
p_{1}(n,1)=p_{2}(n,1)=1
\]
\begin{notation}
	Since the question in question is equivalent to partitions of positive integers , so we use the notation to match the circles on the partitions.\\
	In a circle where there is no circle:  
	\[1\equiv\tikz \draw (0,0) circle [radius=10pt];=B^{1}(1-1)=B^{1}(0) \]
	In general, if there are $n$ circles  in a circle, there are $B(n-1)$ states for that circle.\\
	To get the number of states that the circles are next to each other based on the MP\footnote{The Multiplication Principle(MP)\cite{T}} principle, we multiply them:
	\begin{align*}
		&1+1\equiv\tikz \draw (0,0) circle [radius=10pt];\:\: \tikz \draw (0,0) circle [radius=10pt];=B(0)B(0)=B^{2}(0).
	\end{align*}
	
	Thus, according to the partitions of integers, we can express: 
	\begin{equation}\label{sara 92}
		n =n_1 +n_2 +\cdots+n_l=B(n_1-1)B(n_2-1)\cdots B(n_l-1).
	\end{equation}
\end{notation}
\begin{theorem}\label{sara 79}
	For each number belonging to natural numbers such as $n$, the function $B(n)$ is equal to:
	\[
	B(n)=\sum_{k_{1}+2k_{2}+\cdots+nk_{n}=n}\left(B^{k_{1}}(0)B^{k_{2}}(1)\cdots B^{k_{n}}(n-1)\right).
	\]
\end{theorem}
\begin{proof}
	Given \ref{sara 92} the sum of the states of the partitions is equal to $B(n)$. Therefore, using the equation $ k_{1}+2k_{2}+\cdots+nk_{n}=n $, we can express that:
	\[
	B(n)=\sum_{k_{1}+2k_{2}+\cdots+nk_{n}=n}\left(B^{k_{1}}(1-1)B^{k_{2}}(2-1)\cdots B^{k_{n}}(n-1)\right).
	\]
\end{proof}
Now we try to express another formula for $B(n)$ by obtaining $p_{2}(n,k)$.
\begin{theorem} \label{sara 40}
	 $p_{2}(n,k)$ for every $n$ belongs to natural numbers and for every $k$ belongs to natural numbers where $1<k\le n $: 
	\[
	n=k\left\lfloor \frac{n}{k} \right \rfloor+r \quad,\quad 0\le r \le k-1
	\]
	\[
	p_{2}(n,k)=\left(\sum_{i=1}^{\left \lfloor \frac{n}{k} \right \rfloor-1 }\left(B^{i}(k-1)\sum_{j=1}^{k-1}p_{2}(n-ik,j)\right)\right)+B^{\left \lfloor \frac{n}{k} \right \rfloor}(k-1)B(r)
	\]
\end{theorem}
\begin{proof}
	
	According to the division algorithm, if we assume 
	\[
	\left \lceil \frac{n}
	{2} \right \rceil\le k \le n
	\]
	It can be stated that: 
	\begin{align*}
		n &\equiv B(n-1)B(0)=p_{2}(n,n)\\
		1+(n-1) &\equiv B(n-2)B(1)=p_{2}(n,n-1) \\
		&\vdots\\
		\left(n-\left \lceil \frac{n}
		{2} \right \rceil\right)+\left(\left \lceil \frac{n}
		{2} \right \rceil\right)&\equiv B\left(\left\lfloor \frac{n}{2} \right \rfloor\right)B\left(n-\left \lceil \frac{n}
		{2} \right \rceil\right)=p_{2}\left(n,\left \lceil \frac{n}
		{2} \right \rceil\right)
	\end{align*}
	Hence in the general case for this range $\left \lceil \frac{n}
	{2} \right \rceil\le k \le n$ :
	\[
	p_{2}(n,k)=B(n-k)B(k-1).
	\]
	And now suppose: 
	\[
	2 \le k \le \left\lfloor \frac{n}{2} \right \rfloor
	\]
	In this case, based on the number $k$ , we will have: 
	\begin{align*}
		(n-k)+k &\equiv B^{1}(k-1)\sum_{j=1}^{k-1}p_{2}(n-k,j) \\
		(n-2k)+k+k &\equiv B^{2}(k-1)\sum_{j=1}^{k-1}p_{2}(n-2k,j) \\
		&\vdots\\
		\left(n-\left(\left\lfloor \frac{n}{k} \right \rfloor-1\right)k\right)+k+\cdots+k &\equiv B^{\left\lfloor \frac{n}{k} \right \rfloor-1}(k-1)\sum_{j=1}^{k-1}p_{2}\left(n-\left(\left\lfloor \frac{n}{k} \right \rfloor-1\right)k,j\right) \\
		\left(n-\left\lfloor \frac{n}{k} \right \rfloor k\right)+k+k+\cdots+k &\equiv B^{\left\lfloor \frac{n}{k} \right \rfloor}(k-1)B(r).
	\end{align*}
\end{proof}
\begin{remark}
	According to theorem \ref{sara 40}, it can be stated:\\
	If $ \left\lfloor \frac{n}{k} \right \rfloor=1 $ or $ \left\lceil \frac{n}{2} \right \rceil\le k \le n$ then :
	\[
	p_{2}(n,k)=B(n-k)B(k-1)
	\]
	And if $ \left\lfloor \frac{n}{k} \right \rfloor=2 $  then:
	\[
	p_{2}(n,k)=B(k-1)(p_{2}(n-k,k-1)+\cdots+p_{2}(n-k,1))+B^{2}(k-1)B(r).
	\]	
\end{remark}
\begin{remark} For $k=2$ we will have: 
	\[
	p_{2}(n,2)=\left(\sum_{i=1}^{\left \lfloor \frac{n}{2} \right \rfloor-1 }\left(B^{i}(2-1)\sum_{j=1}^{2-1}p_{2}(n-2i,j)\right)\right)+B^{\left \lfloor \frac{n}{2} \right \rfloor}(2-1)B(r)
	\]
	So we will have: 
	\[
	p_{2}(n,2)=\left\lfloor \frac{n}{2} \right \rfloor.
	\]
\end{remark}
\begin{theorem}
$p_{1}(n,k)$ for every $n$ belongs to natural numbers and for every $k$ belongs to natural numbers where $1<k\le n $:
	\[
	n=k\left\lfloor \frac{n}{k} \right \rfloor+r \quad,\quad 0\le r \le k-1
	\]
	\[p_{1}(n,k)=\left(\sum_{i=1}^{\left\lfloor \frac{n}{k} \right \rfloor-1}\left(\sum_{j=1}^{k-1}p_{1}(n-ik,j)\right)\right)+p(r).\]
\end{theorem}
\begin{proof}
	According to theorem \ref{sara 40}, this theorem is clear.
\end{proof}
\begin{remark}For $k=2$ we will have: 
	\[	p_{1}(n,2)=\left\lfloor \frac{n}{2} \right \rfloor.\]
\end{remark}
\begin{theorem} \label{sara 77}
	For each number belonging to natural numbers such as $n$, the function $B(n)$ is equal to:
	\[
	B(n)=\sum_{k=1}^{n}p_{2}(n,k).
	\]
\end{theorem}
\begin{proof}
	According to theorem \ref{sara 40}, it is clear that  $B(n)$ is equal to the sum of all states $p_{2}(n,k) $. 
\end{proof}
For each $k>n$, $p_{2}(n,k)=p_{1}(n,k)=0$. 
\begin{table}[ht] 
	\caption{$ p_{2}(n,k) $}\label{eqtable}
	\renewcommand\arraystretch{1.5}
	\noindent\[
	\begin{array}{l|cccccccr}
		{k\ n}&{1}&{2}&{3}&{4}&{5}&{6}&{7}&{}\\
		\hline
		{1}&{1}&{1}&{1}&{1}&{1}&{1}&{1}&{}\\
		{2}&{}&{1}&{1}&{2}&{2}&{3}&{3}&{}\\
		{3}&{}&{}&{2}&{2}&{4}&{8}&{10}&{}\\
		{4}&{}&{}&{}&{4}&{4}&{8}&{16}&{}\\
		{5}&{}&{}&{}&{}&{9}&{9}&{18}&{}\\
		{6}&{}&{}&{}&{}&{}&{20}&{20}&{}\\
		{7}&{}&{}&{}&{}&{}&{}&{49}&{}\\
	\end{array}
	\]
\end{table}
\begin{table}[ht] 
	\caption{$ p_{1}(n,k) $}\label{eqtablehh}
	\renewcommand\arraystretch{1.5}
	\noindent\[
	\begin{array}{l|cccccccr}
		{k\ n}&{1}&{2}&{3}&{4}&{5}&{6}&{7}&{}\\
		\hline
		{1}&{1}&{1}&{1}&{1}&{1}&{1}&{1}&{}\\
		{2}&{}&{1}&{1}&{2}&{2}&{3}&{3}&{}\\
		{3}&{}&{}&{1}&{1}&{2}&{3}&{4}&{}\\
		{4}&{}&{}&{}&{1}&{1}&{2}&{3}&{}\\
		{5}&{}&{}&{}&{}&{1}&{1}&{2}&{}\\
		{6}&{}&{}&{}&{}&{}&{1}&{1}&{}\\
		{7}&{}&{}&{}&{}&{}&{}&{1}&{}\\
	\end{array}
	\]
\end{table}
\begin{enumerate}
	\item
\[
p_{1}(n,1)=p_{2}(n,1)=1
\]
\item
\[	p_{1}(n,2)=p_{2}(n,2)=\left\lfloor \frac{n}{2} \right \rfloor\]
\item
\[
p_{1}(n,n-1)=1 \quad,\quad p_{2}(n,n-1)=p_{2}(n-1,n-1)
\]
\item
\[
p_{1}(n,n)=1 \quad,\quad p_{2}(n,n)=B(n-1)=\sum_{k=1}^{n-1}p_{2}(n-1,k)
\]
\item\cite{sara}
\[
p_{1}(n,k)=p_{1}(n-1,k-1)+p_{1}(n-k,k)
\]
\item
\[
p(n)=\sum_{k=1}^{n}p_{1}(n,k) \quad,\quad B(n)=\sum_{k=1}^{n}p_{2}(n,k)\ref{sara 77} 
\]
\item
 $0<x<\frac{1}{4}$
\[
\lim_{n \to \infty} \sum_{k=0}^{n} B(k)x^k=\lim_{n \to \infty} \left(\prod_{k=1}^{n}\frac{1}{1-B(k-1)x^k} \right)
\]
\end{enumerate}
A simple bonuded of this question is obtained using the generating function $(1,S,S^2,...)$, where $S$ is the sum of all possible states of the question. This is the bonuded of the Catalan numbers:
\begin{equation}\label{sara 2}
	B(n)<\frac{\binom{2n}{n}}{n+1}.
\end{equation}
Suppose $x$ is a real number and $0<x<\frac{1}{4}$, in this case according to\eqref{sara 2}: 
\[
B(n)x^{n} < \frac{\binom{2n}{n}}{2^{2n}}\times\frac{1}{n+1}
\]
Given the inequality of\eqref{sara 3}:
\[
B(n)x^{n} < \frac{1}{\sqrt{\pi}}\times\frac{1}{n^{1+\frac{1}{2}}}
\]
Therefore:
\[
\sum_{n=0}^{\infty}B(n)x^{n} < 1+ \sum_{n=1}^{\infty} \left(\frac{1}{\sqrt{\pi}}\times\frac{1}{n^{1+\frac{1}{2}}}\right).
\] 
Series $ \sum_{n=1}^{\infty}\frac{1}{n^{\frac{3}{2}}} $ is a Riemann series that according to the Riemann test is a convergent series, so there is a real number like $w$ that:
\[
\sum_{n=0}^{\infty}B(n)x^{n} < w \quad,\quad w\in\mathbb{R}
\]
By mathematical induction it will be easy to obtain that: 
\[
\forall n\in \mathbb{N} \quad,\quad 2^{n-1} \le B(n)
\]
Therefore, it can be stated that:
\[
\frac{(\sqrt{2})^{2n}}{2} \le B(n) < \frac{2^{2n}}{(n+1)\sqrt{n\pi}}
\]
In which case $ B(n)=(c(n))^{2(n-1)} $ and $ \sqrt{2} \le c(n) < 2 $.

\section{Series}

\begin{theorem} \label{sara 6}
	The generalization of the two series \eqref{sara 189-2} and \eqref{sara 189} is as follows:
	\[
	\forall x \in \mathbb{R}-\{0\}\quad,\quad \forall n\in \mathbb{N}
	\]
	\begin{enumerate}	
		\item
		\[
		\sum_{i=0}^{n} \frac{r}{r+i} (-1)^i \binom{n}{i} x^{r+i}=-\left(\sum_{i=1}^{r}x^{r-i}(1-x)^{n+i}\left(\frac{P_{i}^{r}}{P_{i}^{n+i}}\right)\right)+\frac{1}{\tbinom{n+r}{r}}
		\]
		
		\item
		\[
		\sum_{i=0}^{n} \frac{r}{r+i} \binom{n}{i} x^{r+i}=\left(\sum_{i=1}^{r}(-1)^{i-1}x^{r-i}(1+x)^{n+i}\left(\frac{P_{i}^{r}}{P_{i}^{n+i}}\right)\right)+\frac{(-1)^r}{\binom{n+r}{r}}
		\]
	\end{enumerate}
\end{theorem} 

\begin{proof}[Proof of the first series.\\]

	\[
	S(x,r)=\sum_{i=0}^{n} \frac{r}{r+i} (-1)^i \binom{n}{i} x^{r+i}
	\]
	We obtain the derivative of the series:
	\[
	S^{'}(x,r)=\sum_{i=0}^{n} r (-1)^i \binom{n}{i} x^{r+i-1}
	\]
	We rewrite the phrase in another form:
	\[
	S^{'}(x,r)=rx^{r-1}\sum_{i=0}^{n}\binom{n}{i}  (-1)^i  x^{i}=rx^{r-1}\sum_{i=0}^{n}\binom{n}{i}(-x)^{i}
	\]
	Using a two-sentence series theorem, the phrase is expressed as follows:
	\[
	S^{'}(x,r)=rx^{r-1}(1-x)^{n}
	\]
	Therefore:
	\[
	S(x,r)=\int_{0}^{x} ry^{r-1}(1-y)^{n} dy.
	\]
	Thus, by solving the obtained integral, the correctness of the series can be obtained.\\
	To solve this integral, we use the by parts method repeatedly. Thus it will be obtained:
	\[
	S(x,r)=\left[-\sum_{i=1}^{r-1}y^{r-i}(1-y)^{n+i}\left(\frac{P_{i}^{r}}{P_{i}^{n+i}}\right)-\left(\frac{P_{r}^{r}}{P_{r}^{n+r}}\right)(1-y)^{n+r}\right]_{0}^{x}
	\]
	Now we place the limits of the integral:
	\[
	S(x,r)=-\left(\sum_{i=1}^{r-1}x^{r-i}(1-x)^{n+i}\left(\frac{P_{i}^{r}}{P_{i}^{n+i}}\right)\right)-\frac{(1-x)^{n+r}}{\tbinom{n+r}{r}}+\frac{1}{\tbinom{n+r}{r}}
	\]
	Therefore:
	
	\[
	S(x,r)=-\left(\sum_{i=1}^{r}x^{r-i}(1-x)^{n+i}\left(\frac{P_{i}^{r}}{P_{i}^{n+i}}\right)\right)+\frac{1}{\tbinom{n+r}{r}}.
	\]
\end{proof}
\begin{proof}[	Proof of the second series:]
	\[
	M(x,r)=\sum_{i=0}^{n} \frac{r}{r+i} \binom{n}{i} x^{r+i}
	\]
	We will do the same here as we did to prove the first series.\\
	We obtain the derivative of the series:
	\[
	M^{'}(x,r)=\sum_{i=0}^{n} r \binom{n}{i} x^{r+i-1}
	\]
	We rewrite the phrase in another form: 
	\[
	M^{'}(x,r)=rx^{r-1}\sum_{i=0}^{n} \binom{n}{i} x^{i}=rx^{r-1}(1+x)^{n}
	\]
	Now we express the series using integral:
	\[
	M(x,r)=\int_{0}^{x} ry^{r-1}(1+y)^{n} dy
	\]
	Therefore the result of the integral will be obtained:
	\[
	M(x,r)=\left[\left(\sum_{i=1}^{r-1}(-1)^{i-1}y^{r-i}(1+y)^{n+i}\left(\frac{P_{i}^{r}}{P_{i}^{n+i}}\right)\right)+\frac{(-1)^{r-1}(1+y)^{n+r}}{\binom{n+r}{r}}\right]_{0}^{x}
	\]
	We place the limits of the integral:
	\[
	M(x,r)=\left(\sum_{i=1}^{r}(-1)^{i-1}x^{r-i}(1+x)^{n+i}\left(\frac{P_{i}^{r}}{P_{i}^{n+i}}\right)\right)-(-1)^{r-1}\left(\frac{P_{r}^{r}}{P_{r}^{n+r}}\right)
	\]
	Therefore:
	\[
	M(x,r)=\left(\sum_{i=1}^{r}(-1)^{i-1}x^{r-i}(1+x)^{n+i}\left(\frac{P_{i}^{r}}{P_{i}^{n+i}}\right)\right)+\frac{(-1)^r}{\binom{n+r}{r}}.
	\]	
\end{proof}
\begin{remark}
	The special case of two series in theorem \ref{sara 6}:
	\[
	S(1,r)=\sum_{i=0}^{n}\frac{(-1)^i r}{r+i}\binom{n}{i}=\frac{1}{\tbinom{n+r}{r}}
	\]
	\[
	S(x,1)=-\frac{(1-x)^{n+1}}{n+1}+\frac{1}{n+1}
	\]
	\[
	M(1,r)=\sum_{i=0}^{n} \frac{r}{r+i} \binom{n}{i} = \left(\sum_{i=1}^{r}(-1)^{i-1}2^{n+i}\left(\frac{P_{i}^{r}}{P_{i}^{n+i}}\right)\right)+\frac{(-1)^r}{\binom{n+r}{r}}.
	\]
\end{remark}
\begin{remark}
	Consider the following series:
	\[
	\forall m,n\in\mathbb{N}
	\]
	\[
	\sum_{r=m}^{\infty}\frac{1}{\prod_{i=0}^{n}(r+i)}=\frac{1}{nn!\binom{n+m-1}{m-1}}.
	\]
	The proof of this series is obvious according to the following relation:
	\[
	\frac{n}{r\binom{n+r}{r}}=\frac{1}{\binom{n+r-1}{r-1}}-\frac{1}{\binom{n+r}{r}}.
	\]
\end{remark}

\begin{theorem}
	\[
	\forall m\in \mathbb{R^{+}} \quad,\quad \forall n \in \mathbb{N} \quad,\quad \forall k \in \mathbb{N}
	\]
	\[
	\sum_{r=0}^{n}\binom{n}{r}a^{n-r}r^{k}(bm)^{r}=\sum_{i=1}^{k}S(k,i)P_{i}^{n}(a+bm)^{n-i}(bm)^{i}.
	\]
\end{theorem}
\begin{proof}
	We define the $f(x)$ function:
	\[
	\forall x\in \mathbb{R}\quad,\quad f(x):=(a+be^{x})^n
	\]
	The $ f(x) $ function can be rewritten in another form:
	\[
	f(x)=\sum_{r=0}^{n}\binom{n}{r}a^{n-r}(be^{x})^{r}
	\]
	We get the k-$th$ derivative of the function $ f(x) $ :
	\[
	f^{(k)}(x)=\sum_{r=0}^{n}\binom{n}{r}a^{n-r}r^{k}(be^{x})^{r}
	\]
	Now we get the k-$th$ derivative of the function $ f(x) $  using the general Leibniz rule\cite{S}:
	\[
	f^{(k)}(x)=\sum_{r_{1}+r_{2}+\cdots+r_{n}=k}\binom{k}{r_{1},r_{2},...,r_{n}} (a+be^{x})^{(r_{1})}\cdots (a+be^{x})^{(r_{n})}
	\]
	Now if we consider, we have $ r_{l}\ne 0 $  for $ u $  states, so we have $ r_{l}=0 $ for the remaining $  (n-u)$  states, which can be expressed as: $l,u \in \mathbb{N}_{n}$
	\[
	f^{(k)}(x)=\sum_{u=1}^{k}\left(\sum_{r_{1}+r_{2}+\cdots+r_{n}=k}\binom{k}{r_{1},r_{2},...,r_{n}} (a+be^{x})^{n-u} (be^{x})^{u}\binom{n}{n-u} \right)
	\]
	Using the Stirling numbers of the second kind, it can be stated that:
	\[
	f^{(k)}(x)=S(k,1)1!\binom{n}{1}(a+be^x)^{n-1}be^x+\cdots+S(k,k)k!\binom{n}{k}(a+be^x)^{n-k}b^{k}e^{kx}
	\]
	Therefore:
	\[
	f^{(k)}(x)=\sum_{i=1}^{k}S(k,i)P_{i}^{n}(a+be^{x})^{n-i}(be^{x})^{i}
	\]
	In this case we will have:
	\begin{equation}
		f^{(k)}(x)=\sum_{r=0}^{n}\binom{n}{r}a^{n-r}r^{k}(be^{x})^{r}=\sum_{i=1}^{k}S(k,i)P_{i}^{n}(a+be^{x})^{n-i}(be^{x})^{i}
	\end{equation}
	By limiting the function domain to $\ln(m)$ that $m\in \mathbb{R}^{+}$ , we have:
	\[
	f^{(k)}(\ln(m))=\sum_{r=0}^{n}\binom{n}{r}a^{n-r}r^{k}(bm)^{r}=\sum_{i=1}^{k}S(k,i)P_{i}^{n}(a+bm)^{n-i}(bm)^{i}.
	\]
\end{proof}

\begin{cor}
	Series\eqref{sara 13} is actually the value of $\ln(1)$ in the k-$th$ derivative of the function $f(x)$.
	\[
	f^{(k)}(\ln(1))=\sum_{r=1}^{n}\binom{n}{r} r^k =\sum_{i=1}^{k} S(k,i) P_{i}^{n} 2^{n-i}.
	\]
\end{cor}
\section{Powerful numbers}
\begin{theorem}
	In the following sequence, after $n$ steps of the difference between the sentences of the sequence, we arrive at a fixed sequence of the form $n!$.
	\[
	1^n , 2^n ,3^n ,  ... 
	\]
\end{theorem}
\begin{proof}
	This theorem is equivalent to the number of surjective mappings from $\mathbb{N}_{n}$ to $\mathbb{N}_{n}$ .\\
	The function of the number of surjective  mappings from $\mathbb{N}_{n}$ to $\mathbb{N}_{m}$ is expressed in \cite{T}.
	\[
	F(n,m)=\sum_{k=0}^{m} (-1)^{k}\binom{m}{k} (m-k)^n
	\]
	Since after $n$ steps we get the answer $n!$, so we use the number $n^n$ and the numbers before it. For this reason, we rewrite the sequence in the following form:
	\[
	1^n , 2^n , ... , (n-2)^n ,(n-1)^n , n^n , ...
	\]
	Sentences difference in the first stpe:
	\[
	\dots , n^n-(n-1)^n , (n-1)^n-(n-2)^n , ... , 2^n-1^n 
	\]
	The first step can be shown as follows: 
	\[
	A(n,1)=n^n-(n-1)^n
	\]
	The second step can be expressed:
	\[
	A(n,2)=A(n,1)-A(n-1,1)=n^n-2(n-1)^n-(n-2)^n
	\]
	For the third step:
	\[
	A(n,3)=A(n,2)-A(n-1,2)=n^n-3(n-1)^n-3(n-2)^n-(n-3)^n
	\]
	And the same process can be used for the general case:
	\[
	A(n,n)=A(n,n-1)-A(n-1,n-1)
	\]
	\[
	A(n,n)=\sum_{k=0}^{n}(-1)^k\binom{n}{k}(n-k)^{n}
	\]
	It can be clearly seen that $A(n,n)$ is the same as $F(n,m)$, if $n=m$.
	On the other hand:
	\[
	F(n,m)=m!S(n,m)
	\]
	$ S(n,m) $ is Stirling numbers of the second kind and $ S(n,n)=1  $.\\
	Therefore:
	\[
	A(n,n)=F(n,n)=n!S(n,n)=n!.
	\]
\end{proof}

\bibliographystyle{amsplain}

\end{document}